\documentclass[12pt,reqno]{amsart}

\usepackage[usenames]{color}
\usepackage{amssymb}
\usepackage{amsmath}
\usepackage{amsthm}
\usepackage{amsfonts}
\usepackage{amscd}
\usepackage{graphicx}
\usepackage{mathrsfs}
 
\usepackage{tikz}
\usepackage[colorlinks=true,
linkcolor=webgreen,
filecolor=webbrown,
citecolor=webgreen]{hyperref}

\definecolor{webgreen}{rgb}{0,.5,0}
\definecolor{webbrown}{rgb}{.6,0,0}

\usepackage{color}
\usepackage{fullpage}
\usepackage{float}

\usepackage{graphics}
\usepackage{latexsym}
\usepackage{epsf}
\usepackage{breakurl}

\allowdisplaybreaks
\begin{document}

	\title{	Skew Dyck paths without up--down--left }

	\author[H.~Prodinger]{Helmut Prodinger}
	
	\address{Helmut Prodinger,
		Department of Mathematical Sciences, Stellenbosch University,
		7602 Stellenbosch, South Africa, and
		NITheCS (National Institute for Theoretical and Computational Sciences),
		South Africa}
	\email{hproding@sun.ac.za}

	\date{\today}
	
	\begin{abstract}
Skew Dyck paths without up--down--left are enumerated. In a second step, the number of
contiguous subwords `up--down--left' are counted. This explains and extends results that were
posted in the Encyclopedia of Integer Sequences.
	\end{abstract}
	
	
	\maketitle

	\theoremstyle{plain}
	\newtheorem{theorem}{Theorem}
	\newtheorem{corollary}[theorem]{Corollary}
	\newtheorem{lemma}[theorem]{Lemma}
	\newtheorem{proposition}[theorem]{Proposition}
	
	\theoremstyle{definition}
	\newtheorem{definition}[theorem]{Definition}
	\newtheorem{example}[theorem]{Example}
	\newtheorem{conjecture}[theorem]{Conjecture}
	
	\theoremstyle{remark}
	\newtheorem{remark}[theorem]{Remark}

\section{Introduction}

The entry A128729 in \cite{OEIS} is a bit mysterious. It presents some results about skew Dyck paths without a (contiguous) sequence
`up--down--left,' but without justification, and without reference to a research paper. The present note is designed to fill these gaps.
First, skew Dyck paths are Dyck paths with additional left steps $(-1,-1)$ but without overlapping itself. We found it more convenient~\cite{skew-paper}
to replace a left step by a red down-step $(1,-1)$, labelled red. So we might say that up--down--red is forbidden. In a second step, it won't be forbidden, but counted how often this happens. If it does not happen, it is the same as `forbidden'. Probably the first paper that deals with skew Dyck paths is \cite{Deutsch-italy}.

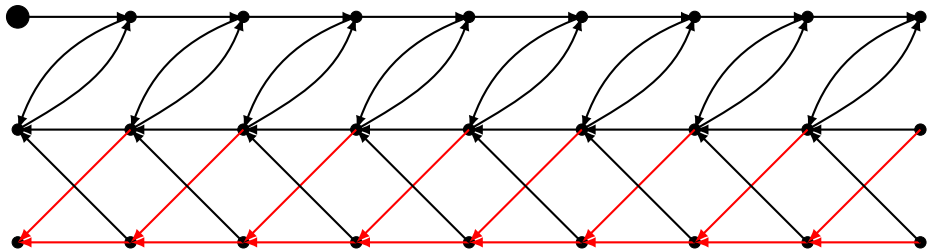
\begin{figure}[h]

		\begin{center}
			\begin{tikzpicture}[scale=1.5]
				\draw (0,0) circle (0.1cm);
				\fill (0,0) circle (0.1cm);
				
				\foreach \x in {0,1,2,3,4,5,6,7,8}
				{
					\draw (\x,0) circle (0.05cm);
					\fill (\x,0) circle (0.05cm);
				}
				
				\foreach \x in {0,1,2,3,4,5,6,7,8}
				{
					\draw (\x,-1) circle (0.05cm);
					\fill (\x,-1) circle (0.05cm);
				}
				
				\foreach \x in {0,1,2,3,4,5,6,7,8}
				{
					\draw (\x,-2) circle (0.05cm);
					\fill (\x,-2) circle (0.05cm);
				}
				
				\foreach \x in {0,1,2,3,4,5,6,7}
				{
					\draw[ thick,-latex] (\x,0) -- (\x+1,0);
					
				}

				\foreach \x in {1,2,3,4,5,6,7}
				{
					\draw[thick,  -latex] (\x+1,0) to[out=200,in=70]  (\x,-1);

				}
				\draw[ thick,     -latex] (1,0) to[out=200,in=70]  (0,-1);

				\foreach \x in {0,1,2,3,4,5,6,7}
				{
					
					\draw[thick,  -latex] (\x,-1) to[out=30,in=250]  (\x+1,0);	
					
				}

				\foreach \x in {0,1,2,3,4,5,6,7}
				{
					\draw[ thick,-latex] (\x+1,-1) -- (\x,-1);
					
				}
				\foreach \x in {0,1,2,3,4,5,6,7}
				{
					\draw[ thick,-latex,red] (\x+1,-1) -- (\x,-2);
					
				}
				
				\foreach \x in {0,1,2,3,4,5,6,7}
				{
					\draw[ thick,-latex,red] (\x+1,-2) -- (\x,-2);
					
				}
				
				\foreach \x in {0,1,2,3,4,5,6,7}
				{
					\draw[ thick,-latex] (\x+1,-2) -- (\x,-1);
					
				}

			\end{tikzpicture}
		\end{center}
		\caption{Three layers of states according to the type of steps leading to them (up, down-black, down-red).}
		\label{arse4}
	\end{figure}

The Figure \ref{arse4} describes skew Dyck paths. One can compute more: Paths that end at level $k$, by any of the 3 types of steps.
Note that up--red and red--up are forbidden, as they would overlap itself. 

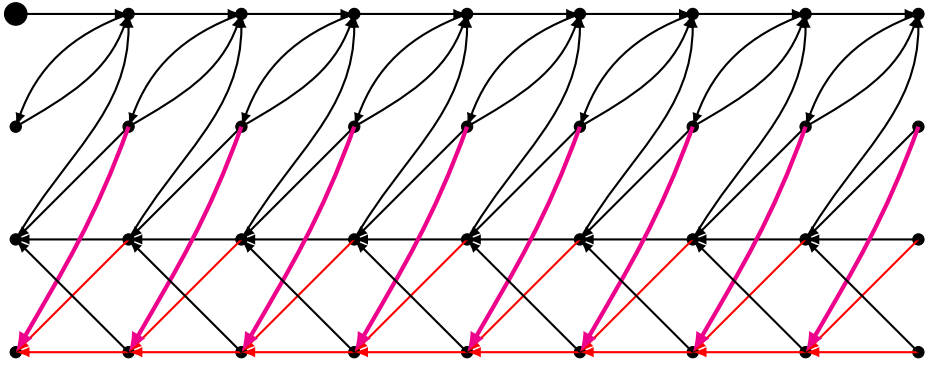
\begin{figure}[h]

	\begin{center}
		\begin{tikzpicture}[scale=1.5]
			\draw (0,0) circle (0.1cm);
			\fill (0,0) circle (0.1cm);
			
			\foreach \x in {0,1,2,3,4,5,6,7,8}
			{
				\draw (\x,0) circle (0.05cm);
				\fill (\x,0) circle (0.05cm);
			}
			
			\foreach \x in {0,1,2,3,4,5,6,7,8}
			{
				\draw (\x,-1) circle (0.05cm);
				\fill (\x,-1) circle (0.05cm);
			}
			
			\foreach \x in {0,1,2,3,4,5,6,7,8}
			{
				\draw (\x,-2) circle (0.05cm);
				\fill (\x,-2) circle (0.05cm);
			}
			\foreach \x in {0,1,2,3,4,5,6,7,8}
			{
				\draw (\x,-3) circle (0.05cm);
				\fill (\x,-3) circle (0.05cm);
			}

			\foreach \x in {0,1,2,3,4,5,6,7}
			{
				\draw[ thick,-latex] (\x,0) -- (\x+1,0);
				
			}

			\foreach \x in {1,2,3,4,5,6,7}
			{
				\draw[thick,  -latex] (\x+1,0) to[out=200,in=70]  (\x,-1);

			}
			\draw [thick,      -latex] (1,0) to[out=200,in=70]  (0,-1);

			\foreach \x in {0,1,2,3,4,5,6,7}
			{
				
				\draw[thick,  -latex] (\x,-1) to[out=30,in=250]  (\x+1,0);	
				
			}

			\foreach \x in {0,1,2,3,4,5,6,7}
			{
				\draw[ thick,-latex,red] (\x+1,-2) -- (\x,-3);
				
			}
			
			\foreach \x in {0,1,2,3,4,5,6,7}
			{
				\draw[ thick,-latex,red] (\x+1,-3) -- (\x,-3);
				
			}
			
			\foreach \x in {0,1,2,3,4,5,6,7}
			{
				\draw[ thick,-latex] (\x+1,-2) -- (\x,-2);
				
			}

			\foreach \x in {0,1,2,3,4,5,6,7}
			{
				\draw[ thick,-latex] (\x+1,-1) -- (\x,-2);
				
			}
			\foreach \x in {0,1,2,3,4,5,6,7}
			{
				\draw[ thick,-latex] (\x,-2) to [out=60,in=-90] (\x+1,-0);
				
			}
			
				\foreach \x in {0,1,2,3,4,5,6,7}
			{
				\draw[ ultra thick, magenta ,-latex] (\x+1,-1) to [out=-110,in=60] (\x,-3);
				
			}
		\foreach \x in {0,1,2,3,4,5,6,7}
		{
			\draw[   thick ,-latex] (\x+1,-3) to (\x,-2);
			
		}

		\end{tikzpicture}
	\end{center}
	\caption{Four layers of states according to the type of steps leading to them.}
	\label{schoas}
\end{figure}

The next step is to take care of up--down--red. Such a red step is depicted in Figure \ref{schoas} in a special color.

\section{The magenta colored steps are forbidden}

 We introduce 4 sequences of generating functions, according to the 4 layers in the defining graph. The variable $z$ counts
 the steps. The recursions can be seen immediately from the graph: 
\begin{align*}
f_{n+1}&=zf_n+zg_n+zh_n,\ f_0=1,\\
g_{n}&=zf_{n+1},\\
h_n&=zg_{n+1}+zh_{n+1}+zk_{n+1},\\
k_n&=zh_{n+1}+zk_{n+1}.
\end{align*}
And now we use a second variable $u$ to be able to work with bivariate generating functions;
 \begin{align*}
 	\sum_{n\ge0}u^{n+1}f_{n+1}&=\sum_{n\ge0}zu^{n+1}f_n+\sum_{n\ge0}zu^{n+1}g_n+\sum_{n\ge0}zu^{n+1}h_n, \\
 	\sum_{n\ge0}u^{n+1}g_{n} &=\sum_{n\ge0}zu^{n+1}f_{n+1},\\
 	\sum_{n\ge0}u^{n+1}h_n&=\sum_{n\ge0}zu^{n+1}g_{n+1}+\sum_{n\ge0}zu^{n+1}h_{n+1}+\sum_{n\ge0}zu^{n+1}k_{n+1},\\
 	\sum_{n\ge0}u^{n+1}k_n&=\sum_{n\ge0}zu^{n+1}h_{n+1}+\sum_{n\ge0}zu^{n+1}k_{n+1}.
 \end{align*}
Translating this into bivariate generating functions $F(u)=F(u,z)$ etc., leads to
\begin{align*}
	F(u)&=1+zuF(u) +zuG(u) +zuH(u), \\
	uG(u)  &=zF(u)-z, \\
	uH(u) &=zG(u)-zg_0 +zH(u)-zh_0 + zK(u)-zk_0 ,  \\
	uK(u) &=zH(u)-zh_0 + zK(u)-zk_0.
\end{align*}

Solving the system leads to
\begin{align*}
F(u)&=\frac{{u}^{2}{z}^{2}g_0+{u}^{2}{z}^{2}h_0+{u}^{2}{z}^{2}k_0-u
	{z}^{3}+2 zu-{u}^{2}+{u}^{2}{z}^{2}-u{z}^{3}g_0-{z}^{4}
}{-{u}^{2}+z{u}^{3}+2zu-{u}^{2}{z}^{2}-u{z}^{3}-{z}^{4}},\\
G(u)&=\frac{ {z}^{2} \left( uzg_0+uzh_0+uzk_0+2 zu-{z}^{2}g_0-
	{u}^{2} \right)
}{-{u}^{2}+z{u}^{3}+2zu-{u}^{2}{z}^{2}-u{z}^{3}-{z}^{4}},\\
H(u)&=\frac{ -z \left( u{z}^{2}-{z}^{3}-g_0 u-h_0 u-k_0 u+z g_0+z{u}^{2}g_0+{u}^{2}zh_0+{u}^{2}zk_0+u{z}^{2}h_0+u{z}^{2}k_0-{z}^{3}g_0 \right)
}{-{u}^{2}+z{u}^{3}+2zu-{u}^{2}{z}^{2}-u{z}^{3}-{z}^{4}},\\
K(u)&=\frac{- \left( {u}^{2}zk_0+{u}^{2}zh_0-h_0 u+u{z}^{2}g_0-k_0 u+u{z}^{2}k_0+u{z}^{2}h_0-zg_0+{z}^{3}+{z}^
	{3}g_0+k_0 {z}^{3}+h_0 {z}^{3} \right) z
}{-{u}^{2}+z{u}^{3}+2zu-{u}^{2}{z}^{2}-u{z}^{3}-{z}^{4}},\\
\end{align*}
One cannot just plug in $u=0$ to identify the constants. However,
\begin{equation*}
-{u}^{2}+z{u}^{3}+2zu-{u}^{2}{z}^{2}-u{z}^{3}-{z}^{4}=z(u-u_1)(u-u_2)(u-u_3),
\end{equation*}
and the bad factors $(u-u_2)(u-u_3)$ can be cancelled out, both, in the denominator and the numerator. 
This is essential for the kernel method \cite{prodinger-kernel}. This leads to
\begin{align*}
	F(u)&=\frac{z^2g_0+z^2h_0+z^2k_0-1+z^2}{z(u-u_1)},\\
	G(u)&=\frac{-z }{(u-u_1)},\\
	H(u)&=\frac{-z(g_0+h_0+k_0)}{(u-u_1)},\\
	K(u)&=\frac{-z(h_0+k_0)}{(u-u_1)}.
\end{align*}
Now we can plug in $u=0$ and solve to get
\begin{align*}
g_0&=\frac{z}{u_1},\\
h_0&=\frac{1-z^2-zu_1}{zu_1},\\
k_0&=\frac{1-z^2-zu_1}{u_1(u_1-z)}.
\end{align*}
In total
\begin{align*}
1+g_0+h_0+k_0&=\frac{ u_1-z^2u_1-z^3}{zu_1(u_1-z)}=\frac{1-zu_1}{z^2}\\
&=1+z^2+2 z^4+6 z^6+20 z^8+71 z^{10}+262 z^{12}+994 z^{14}+3852 z^{16}+\cdots.
\end{align*}
The last simplification was done by noticing that $u_1=\textsf{RootOf}( 
-u^2+zu^3+2zu-u^2z^2-uz^3-z^4,u)$
and using Maple's \textsf{evala} command. (One can also use the \textsf{algeqtoseries} command to go from the algebraic equation to the series expansion.)
  
 We further find
 \begin{equation*}
F+G+H+K=-\frac{u_1(1-zu_1)}{z^2(u-u_1)}
 \end{equation*}
and
 \begin{equation*}
	[u^k](F+G+H+K)=\frac{1-zu_1}{z^2u_1^{k}}.
\end{equation*}

We summarize the results of this section:
\begin{theorem}
	The generating function of skew Dyck paths without up--down--left ending on level $k$ is given by
	\begin{equation*}
		 \frac{1-zu_1}{z^2u_1^{k}}
	\end{equation*}
where
\begin{equation*}
u_1=\frac1z -z-z^3-2z^5-6z^7-20z^9-71z^{11}-262z^{13}-\cdots
\end{equation*}
is the solution of the algebraic equation $-u^2+zu^3+2zu-u^2z^2-uz^3-z^4=0$ that does not contain imaginary numbers in its series expansion around $z=0$. Equivalently, the other solutions $u_2$, $u_3$ are such that $1/(u-u_2)$ and $1/(u-u_3)$ have no power series expansion around $(0,0)$.
\end{theorem}

Noticing that one can only return to level 0 in an even number of steps, one can write an alternative form for
$f_0+g_0+h_0+k_0$:
\begin{align*}
	\textsf{RootOf}&(z^2u^3-z(2-z)u^2+(1-z^2)u-1+z+z^2 ,u)\\
	&=1+z+2 z^2+6 z^3+20 z^4+71 z^{5}+262 z^{6}+994 z^{7}+3852 z^{8}+\cdots.
\end{align*}
We show how the algebraic equation for $u_1$ is transformed into one for $\frac{1-zu_1}{z^2}$; in the last step, $z^2=Z$ was substituted:
\begin{gather*}
u_1 \longrightarrow -u^2+zu^3+2zu-u^2z^2-uz^3-z^4=0,\\
zu_1 \longrightarrow -(zu)^2+(zu)^3+2z^2(zu)-z^2(zu)^2-z^4(zu)-z^6=0,\\
-zu_1 \longrightarrow -(-zu)^2-(-zu)^3-2z^2(-zu)-z^2(-zu)^2+z^4(-zu)-z^6=0,\\
1-zu_1 \longrightarrow -(1-zu-1)^2-(1-zu-1)^3-2z^2(1-zu-1)\\
\qquad\qquad\qquad\qquad\qquad-z^2(1-zu-1)^2+z^4(1-zu-1)-z^6=0,\\
\qquad = 2(1-zu)^2-(1-zu)-(1-zu)^3+z^2-z^2(1-zu)^2+z^4(1-zu)-z^4-z^6=0,\\
\frac{1-zu_1}{z^2}=U_1\longrightarrow
2z^4U^2-z^2U-z^6U^3+z^2-z^6U^2+z^6U-z^4-z^6=0,\\
U_1\longrightarrow
2z^2U^2-U-z^4U^3+1-z^4U^2+z^4U-z^2-z^4=0,\\
U_1\longrightarrow
2ZU^2-U-Z^2U^3+1-Z^2U^2+Z^2U-Z-Z^2=0.
\end{gather*}

So, writing $S(z)=f_0+g_0+h_0+k_0$, the algebraic equation is
\begin{equation*}
z^2S^3 - z(2-z)S^2 + (1 - z^2)S - 1 + z + z^2 = 0,
\end{equation*}
this is the form given in A128729 \cite{OEIS}, and $z$ counts only the half-length of the skew Dyck path,
and Maple's gfun translates that into a differential equation
\begin{equation*}
31z-8-15zS-(2z-1)(44z^3+15z^2-48z+8)S'-z(11z^2+16z-4)(2z-1)^2S{''}=0,
\end{equation*}
from which gfun derives a linear recursion with polynomial coefficients:
\begin{multline*}
-44n(n+1)s_n-2(n+1)(10n-7)s_{ n+1}\\+3(115+106n+23n^2)s_{n+2}-32(n+4)(n+3)s_{n+3}+4(n+5)(n+4)s_{n+4}=0,
\end{multline*}
where
\begin{equation*}
S(z)=\sum_{n\ge0}s_nz^n.
\end{equation*}

\subsection*{Asymptotics}

The technique is \emph{singularity analysis of generating functions}, as described in \cite{FlOd90}.

We consider $eq:=z^2S^3-z(2-z)S^2+(1-z^2)S-1+z+z^2 $ and find the dominant singularity $(z_0,S_0)$.
For that, we consider the equation 
\begin{equation*}
\frac{eq}{dS}=3z^2S^2-2z(2-z)S+1-z^2=0
\end{equation*}
and thus $S_0=\frac{z+1}{3z}$. Plugging this into $eq$, we find $z_0=\frac{2}{11}(3\sqrt3-4)$ as closest solution to the origin. From this, $S_0=1+\frac{\sqrt3}{2}$. Expanding $eq$ locally,
\begin{equation*}
z-z_0\sim \frac{216-129\sqrt3}{121}\Big(S-1-\frac{\sqrt3}{2}\Big)^2,
\end{equation*}
or
\begin{equation*}
	(z_0-z)\Big(8+\frac{43}{9}\sqrt3\Big) =\Big(2+\frac89\sqrt3\Big)\Big(1-\frac{z}{z_0}\Big)\sim  \Big(S-1-\frac{\sqrt3}{2}\Big)^2.
\end{equation*}
This gives us the local expansion around the square-root singularity:
\begin{equation*}
S\sim 1+\frac{\sqrt3}{2}-\sqrt{2+\frac 89\sqrt3}\;\sqrt{1-\frac {z}{z_0}}.
\end{equation*}
Therefore
\begin{equation*}
[z^n]S\sim \sqrt{2+\frac 89\sqrt3}\;\frac{1}{2\sqrt\pi}\Big(2+\frac32\sqrt3\Big)^{n}n^{-3/2}.
\end{equation*}
\section{Counting the up--down--red configurations}

This is not too different from what we did before, but now we use another variable, $t$, which is attached to the magenta colored edges. In this way, the exponent of $t$ 
refers to the number of up--down--red configurations;
\begin{align*}
	f_{n+1}&=zf_n+zg_n+zh_n,\ f_0=1,\\
	g_{n}&=zf_{n+1},\\
	h_n&=zg_{n+1}+zh_{n+1}+zk_{n+1},\\
	k_n&=ztg_{n+1}+zh_{n+1}+zk_{n+1}.
\end{align*}
We use again a second variable $u$ to be able to work with bivariate generating functions;
\begin{align*}
	\sum_{n\ge0}u^{n+1}f_{n+1}&=\sum_{n\ge0}zu^{n+1}f_n+\sum_{n\ge0}zu^{n+1}g_n+\sum_{n\ge0}zu^{n+1}h_n, \\
	\sum_{n\ge0}u^{n+1}g_{n} &=\sum_{n\ge0}zu^{n+1}f_{n+1},\\
	\sum_{n\ge0}u^{n+1}h_n&=\sum_{n\ge0}zu^{n+1}g_{n+1}+\sum_{n\ge0}zu^{n+1}h_{n+1}+\sum_{n\ge0}zu^{n+1}k_{n+1},\\
	\sum_{n\ge0}u^{n+1}k_n&=\sum_{n\ge0}ztu^{n+1}g_{n+1}+\sum_{n\ge0}zu^{n+1}h_{n+1}+\sum_{n\ge0}zu^{n+1}k_{n+1}.
\end{align*}
Translating this into bivariate generating functions $F(u)=F(u,z)$ etc., leads to
\begin{align*}
	F(u)&=1+zuF(u) +zuG(u) +zuH(u), \\
	uG(u)  &=zF(u)-z, \\
	uH(u) &=zG(u)-zg_0 +zH(u)-zh_0 + zK(u)-zk_0 ,  \\
	uK(u) &=ztG(u)-ztg_0 +zH(u)-zh_0 + zK(u)-zk_0.
\end{align*}
Solving the system leads to
\begin{align*}
F(u)&=\frac{\mathscr{F}}{2zu-z^4-u^2+z^4t-uz^3+zu^3-u^2z^2},\\
G(u)&=\frac{\mathscr{G}}{2zu-z^4-u^2+z^4t-uz^3+zu^3-u^2z^2},\\
H(u)&=\frac{\mathscr{H}}{2zu-z^4-u^2+z^4t-uz^3+zu^3-u^2z^2},\\
K(u)&=\frac{\mathscr{K}}{2zu-z^4-u^2+z^4t-uz^3+zu^3-u^2z^2},
 \end{align*}
with
$\mathscr{F}={u}^{2}{z}^{2}h_0+{z}^{3}utg_0+2 zu-{u}^{2}+{z}^{4}t-{z}^{4}+{u}^{2}{z}^{2}k_0
+{u}^{2}{z}^{2}-u{z}^{3}+{u}^{2}{z}^{2}g_0-u{z}^{3}g_0$,
$\mathscr{G}={z}^{3}uh_0+{z}^{4}g_0 t+{z}^{3}uk_0+2 u{z}^{3}+u{z}^
{3}g_0-{z}^{4}g_0-{u}^{2}{z}^{2}$,
$\mathscr{H}={z}^{2}tg_0-{z}^{3}utg_0-{z}^{4}g_0 t+{z}^{4}-u{z}^{3}
-{u}^{2}{z}^{2}h_0-{z}^{4}t+uzg_0+uzh_0+uzk_0-{z}^
{2}g_0-{u}^{2}{z}^{2}k_0-{u}^{2}{z}^{2}g_0-{z}^{3}uh_0-{z}^{3}uk_0+{z}^{4}g_0$,
$\mathscr{K}=uztg_0-{z}^{2}{u}^{2}tg_0-{z}^{2}tg_0+{z}^{4}g_0 
t+{z}^{4}h_0 t+{z}^{4}k_0 t-{z}^{4}-{z}^{3}ut+uzh_0+u
zk_0+{z}^{2}g_0-{z}^{4}g_0-u{z}^{3}g_0-{z}^{4}h_0-{z}^{3}uh_0-{u}^{2}{z}^{2}h_0-{z}^{4}k_0-{z}^{3
}uk_0-{u}^{2}{z}^{2}k_0+{z}^{4}t
$.

The denominator $2zu-z^4-u^2+z^4t-uz^3+zu^3-u^2z^2$ still factors as $z(u-u_1)(u-u_2)(u-u_3)$ where $u_1$, $u_2$, $u_3$ now depend on $t$.
Again, the factor $(u-u_2)(u-u_3)$ can be cancelled out:
\begin{align*}
	F(u)&=\frac{z^2g_0+z^2h_0+z^2k_0-1+z^2}{z(u-u_1)},\\
	G(u)&=\frac{-z }{(u-u_1)},\\
	H(u)&=\frac{-z(g_0+h_0+k_0)}{(u-u_1)},\\
	K(u)&=\frac{-z(tg_0+h_0+k_0)}{(u-u_1)}.
\end{align*}
Now we can plug in $u=0$ and compute the constants:
\begin{align*}
g_0&=\frac{z}{u_1},\\
h_0&=\frac{1-z^2-zu_1}{zu_1},\\
k_0&=\frac{(1-z^2-zu_1)(tu_1-tz+z)}{zu_1(u_1+tz-z)}.
\end{align*}
In total
\begin{equation*}
1+g_0+h_0+k_0=\frac{-tzu_1^2+tu_1+z^3t-z^2u_1-z^3+u_1}{zu_1(u_1+tz-z)}=\frac{1-zu_1}{z^2}.
\end{equation*}
The series starts as
\begin{equation*}
1+z^2+(2+t)z^4+(4t+6)z^6+(16t+20)z^8+(71+64t+2t^2)z^{10}+(262+261t+20t^2)z^{12}+\cdots.
\end{equation*}
One can compute the triple generating function, where $z$ refers to the length of the walk, $u$ to the level at the end, and $t$ to the number
of up--down--red configurations:
\begin{equation*}
F(u)+G(u)+H(u)+K(u)=\frac{z^4t-u_1^2+tz^2u_1^2+z^3u_1+z^2u_1^2-2u_1tz-t^2z^4}{u_1(u_1+tz-z)z(u-u_1)}=\frac{u_1(zu_1-1)}{z^2(u-u_1)}.
\end{equation*}
Consequently
\begin{equation*}
[u^k]	(F(u)+G(u)+H(u)+K(u))=\frac{1-zu_1}{z^2u_1^{k}}.
\end{equation*}
Since only an even number of steps brings us to level 0, we have
\begin{equation*}
z^2R^3 - z(2-z)R^2 + (1-z^2)R - 1 + z + z^2 - tz^2 = 0
\end{equation*}
for $1+g_0+h_0+k_0$ where $z$ counts the half-length of the skew Dyck path. This is the form given in  sequence A128728~\cite{OEIS}.

\end{document}